\documentclass[reqno]{amsart}
\usepackage{hyperref}
\newtheorem*{theoA}{Theorem A}
\newtheorem*{theoB}{Theorem B}
\newtheorem*{theoC}{Theorem C}
\newtheorem*{theoD}{Theorem D}
\newtheorem*{theoE}{Theorem E}

\newtheorem{theo}{Theorem}[section]
\newtheorem{lem}{Lemma}[section]

\newtheorem{ques}{Question}[section]
\newtheorem{exm}{Example}[section]
\newtheorem{defi}{Definition}[section]
\newtheorem{rem}{Remark}[section]
\newcommand{\ol}{\overline}
\newcommand{\be}{\begin{equation}}
\newcommand{\ee}{\end{equation}}
\newcommand{\beas}{\begin{eqnarray*}}
	\newcommand{\eeas}{\end{eqnarray*}}
\newcommand{\bea}{\begin{eqnarray}}
\newcommand{\eea}{\end{eqnarray}}

\numberwithin{equation}{section}

\begin{document}
\title[\hfilneg \hfil  Uniqueness of q-differential polynomials]
{ Value distribution and uniqueness for q-difference of meromorphic functions Sharing Two Sets}
\author[Goutam Haldar\hfil
\hfilneg]
{Goutam Haldar}



\subjclass[2010]{30D35.}
\keywords{ Meromorphic function, q-difference operator, small function, weighted sharing, zero-order}

\maketitle

\begin{abstract}
   In this paper, we investigate the value distribution for linear q-difference polynomials of transcendental meromorphic functions of zero order which improves the results of Xu, Liu and Cao (\cite{Xu & Liu & Cao & 2015}). We also investigate the uniqueness of zero order meromorphic function with its q-difference operator sharing two sets with finite weight. Some examples have been exhibited which are relevant to the content of the paper.
   \end{abstract}

\section{\textbf{Introduction}}

Let $f$ and $g$ be two non-constant meromorphic functions defined in the open complex plane $\mathbb{C}$. If for some $a\in\mathbb{C}\cup\{\infty\}$, the zero of  $f-a$ and $g-a$ have the same locations as well as same multiplicities, we say that $f$ and $g$ share the value $a$ CM (counting multiplicities). If we do not consider the multiplicities, then $f$ and $g$ are said to share the value $a$ IM (ignoring multiplicities). We adopt the standard notations of the Nevanlinna theory of meromorphic functions explained in (\cite{Hayman & 1964}).\par A meromorphic function $a$ is said to be a small function of $f$ provided that $T(r,a)=S(r,f)$. i.e., $T(r,a)=o(T(r,f))$ as $r\longrightarrow \infty$, outside of a possible exceptional set of finite linear measure.\par
For a set $S\subset\mathbb{C}$, we define
\begin{equation*} E_{f}(S)=\bigcup_{a\in S}\{z|f(z)=a(z)\},\end{equation*} where each zero is counted according to its multiplicity and $$\overline E_{f}(S)=\bigcup_{a\in S}\{z|f(z)=a(z)
\},\;\text{where each zero is counted only once}.$$\par
If $E_{f}(S)=E_{g}(S)$, we say that $f$, $g$ share the set S CM and if $\overline E_{f}(S)=\overline E_{g}(S)$, we say $f$, $g$ share the set S IM.\par
In 2001, Lahiri (\cite{Lahiri & Nagoya & 2001}) introduced a gradation of sharing of values or sets which is known as weighted sharing. Below we are recalling the notion.
\begin{defi}(\cite{Lahiri & Nagoya & 2001})
	Let $k$ be a non-negative integer or infinity. For $a\in \mathbb{C}\cup\{\infty\}$ we denote by $E_{k}(a,f)$ the set of all $a$-points of $f$, where an $a$ point of multiplicity $m$ is counted $m$ times if $m\leq k$ and $k+1$ times if $m>k.$ If $E_{k}(a,f)=E_{k}(a,g),$ we say that $f$, $g$ share the value $a$ with weight $k$. \end{defi}
We write $f$, $g$ share $(a,k)$ to mean that $f,$ $g$ share the value $a$ with weight $k.$ Clearly if $f,$ $g$ share $(a,k)$ then $f,$ $g$ share $(a,p)$ for any integer $p$, $0\leq p<k.$ Also we note that $f,$ $g$ share a value $a$ IM or CM if and only if $f,$ $g$ share $(a,0)$ or $(a,\infty)$ respectively.\par
\begin{defi}\cite{Lahiri & Sarkar & 2004}Let $p$ be a positive integer and $a\in\mathbb{C}\cup\{\infty\}$.\begin{enumerate}
		\item[(i)] $N(r,a;f\mid \geq p)$ ($\ol N(r,a;f\mid \geq p)$)denotes the counting function (reduced counting function) of those $a$-points of $f$ whose multiplicities are not less than $p$.\item[(ii)]$N(r,a;f\mid \leq p)$ ($\ol N(r,a;f\mid \leq p)$)denotes the counting function (reduced counting function) of those $a$-points of $f$ whose multiplicities are not greater than $p$.
	\end{enumerate}
\end{defi}
\begin{defi} \cite {Alzahare & Yi & 2004} Let $f$ and $g$ be two non-constant meromorphic functions such that $f$ and $g$ share the value $a$ IM. Let $z_{0}$ be a $a$-point of $f$ with multiplicity $p$, a $a$-point of $g$ with multiplicity $q$. We denote by $\ol N_{L}(r,a;f)$ the counting function of those $a$-points of $f$ and $g$ where $p>q$, by $N^{1)}_{E}(r,a;f)$ the counting function of those $a$-points of $f$ and $g$ where $p=q=1$ and by $\ol N^{(2}_{E}(r,a;f)$ the counting function of those $a$-points of $f$ and $g$ where $p=q\geq 2$, each point in these counting functions is counted only once. Similarly, one can define $\ol N_{L}(r,a;g),\; N^{1)}_{E}(r,a;g),\; \ol N^{(2}_{E}(r,a;g).$
\end{defi}
\begin{defi}\cite{Lahiri & Nagoya & 2001,Lahiri & Complex Var & 2001} Let $f$, $g$ share a value $a$ IM. We denote by $\ol N_{*}(r,a;f,g)$ the reduced counting function of those $a$-points of $f$ whose multiplicities differ from the multiplicities of the corresponding $a$-points of $g$.
	Clearly $\ol N_{*}(r,a;f,g)\equiv\ol N_{*}(r,a;g,f)$ and $\ol N_{*}(r,a;f,g)=\ol N_{L}(r,a;f)+\ol N_{L}(r,a;g)$.
\end{defi}
Recently, many of authors have shown their interest in studying difference equations, the difference product and the q-difference analogues the value distribution theory in the complex plane $\mathbb{C}$. A number of remarkable research works (see \cite{Chen & Huang & Zheng & 2011}, \cite{Chiang & Feng & 2008}, \cite{Halburd & Korhonen & 2006}, \cite{Halburd & Korhonen & Acad.2006}, \cite{Heittokangas et al & 2011}, \cite{Laine & Yang & 2007}, \cite{Liu & 2009}, \cite{Liu & Yang & 2009}, \cite{Liu & Liu & Cao & 2011}, \cite{Luo & Lin &  2011}, \cite{Wang & Xu & Zhang & 2014}, \cite{Zhang & 2010}) have focused on the uniqueness of difference analogues of Nevanlinna theory. In 2006, Halburd and Korhonen (\cite{Halburd & Korhonen & 2006}) established a difference analogue of the Logarithmic Derivative Lemma, and then applying it, a lot of results on meromorphic solutions of complex difference equations has been proved. After that Barnett, Halburd, Korhonen and Morgan (\cite{Barnett & Halburd & 2007}) also established a q-difference analogue of the Logarithmic Derivative Lemma. \par 
Let us first recall the notion of the q-shift and q-difference operator of a meromorphic function $f$.
\begin{defi}
	For a meromorphic function f and $c,\; q(\neq 0)\in\mathbb{C}$, let us now denote its q-shift $E_{q,c}f$ and q-difference operators $\Delta_{q,c}f$ respectively by $E_{q,c}f(z) = f(qz+c)$ and $\Delta_{q,c}f(z)=f(qz+c)-f(z)$, $\Delta_{q,c}^{k}f(z) := \Delta_{q,c}^{k-1}(\Delta_{q,c}f(z)),\; \text{for all}\;  k\in\mathbb{N}-\{1\}$.  \end{defi}   
We now proceed to define linear q-shift and q-difference operators, denoted respectively by $L_{k}(f, E_q)$ and $L_{k}(f,\Delta)$ of a meromorphic function $f$ in a more compact and convenient way in the following, which is one of the motivation of writing the paper.
\begin{defi}
	Let us define
\begin{equation} \label{e1.1a}L_{k}(f, E_q)=a_kf(q_kz+c_k)+a_{k-1}f(q_{k-1}z+c_{k-1})+\ldots+a_0f(q_0z+c_0)\end{equation} and	
\begin{equation} \label{e1.1}L_{k}(f, \Delta)=a_{k}\Delta_{q_k,c_k}f(z)+a_{k-1}\Delta_{q_{k-1},c_{k-1}}f(z)+\ldots+a_{0}\Delta_{q_{0},c_{0}}f(z),\end{equation} where $a_0, a_1, \ldots, a_{k}$; $q_0, q_1, \ldots, q_k; c_0,c_1,\ldots,c_k$ are complex constants.
From (\ref{e1.1}), one can easily observe that \begin{equation*} L_{k}(f, \Delta)= L_{k}(f, E_q)-\sum_{j=0}^{k}a_j f(z).\end{equation*}\end{defi}
If we choose $q_{j}=q^{j}$, $c_j=c$ and $a_{j}=(-1)^{k-j}{k\choose j}$ for $0\leq j\leq k$, then $L_{k}(f, \Delta)$ reduces $\Delta_{q,c}^{k}f(z)$. \par 
Let $P(z)=a_nz^n+a_{n-1}z^{n-1}+\ldots+a_0$ be a nonzero polynomial of degree $n$, where $a_{n}(\neq 0), a_{n-1}, \ldots, a_0$ are complex constants and $m$ be the number of distinct zeros of $P(z)$. From now on,  unless otherwise stated, for $a\neq0$, we denote, by $S_1=\{a, aw, aw^2,\ldots, aw^{n-1}\}$, where $n\in\mathbb{N}$, $w$ is a n-th roots of unity and $S_2=\{\infty\}$ throughout the paper.\par 
Zhang and Korhonen  (\cite{Zhang & Korhonen & 2010}) studied the value distribution of q-difference polynomials of meromorphic functions and obtained the following result.
\begin{theoA}\cite{Zhang & Korhonen & 2010}
	Let $f$ be a transcendental meromorphic (resp. entire) function of zero order and q non-zero complex constant. Then for $n\geq6$	(resp. $n \geq2$), $f(z)^nf(qz)$ assumes every non-zero value $a\in\mathbb{C}$ infinitely often.\end{theoA}
Recently, Liu and Qi \cite{Liu & Qi & 2011} firstly investigated value distributions for a q-shift of the meromorphic function and obtained the following result.
\begin{theoB}\cite{Liu & Qi & 2011}
	Let $f$ be a zero-order transcendental meromorphic function, $n\geq6$; $q\in\mathbb{C}-\{0\}$, $\eta\in\mathbb{C}$ and $R(z)$ a rational function. Then $f(z)^nf(qz + \eta)-R(z)$ has infinitely many zeros.\end{theoB}
In 2015, Xu, Liu and Cao \cite{Xu & Liu & Cao & 2015} started investigation about the zeros of $P(f)f(qz + \eta) = a(z)$and $P(f)[f(qz + \eta)-f(z)] = a(z)$, where $a(z)$ is a small function of $f$ and obtained the following two results. \par 
\begin{theoC}\cite{Xu & Liu & Cao & 2015}
	Let $f$ be a zero-order transcendental meromorphic (resp. entire) function, $q(\neq 0), \eta$ are complex constants. Then for $n > m+4$ (resp. $n > m$), $P(f)f(qz+\eta)=a(z)$ has infinitely many solutions, where a(z) is a non-zero small functions in $f$.\end{theoC}
\begin{theoD}\cite{Xu & Liu & Cao & 2015}
	Let $f$ be a zero-order transcendental meromorphic (resp. entire) function, $q(\neq 0), \eta$ are complex constants. Then for $n > m+6$ (resp. $n > m+2$), $P(f)\{f(qz+\eta)-f(z)\}=a(z)$ has infinitely many solutions, where a(z) is a non-zero small functions in $f$.\end{theoD}

Regarding uniqueness, in 2011, Qi, Liu and Yang \cite{Qi & Liu & Yang & 2011} obtained the following result.
\begin{theoE}\cite{Qi & Liu & Yang & 2011}
	Let $f$ be a zero-order meromorphic function, and $q\in\mathbb{C}-\{0\}$, $n\geq4$ be an integer, and let $F=f^n$. If $F(z)$ and $F(qz)$ share $a\in\mathbb{C}-\{0\}$ and	$\infty$ CM, then $f(z) = tf (qz)$ for a constant $t$ that satisfies $t^n=1$.
\end{theoE}
\begin{ques}\label{q1.1a}
	What would happen if we replace $P(f)f(qz+c)$ in Theorem C by more general $q$-shift difference polynomial $P(f)L_k(f,E_k)$ and $P(f)\{f(qz+c)-f(z)\}$ in Theorem D by $P(f)L_k(f,\Delta)$?
\end{ques}
\begin{ques}\label{q1.1}
	what would happen if we replace $F=f^{n}$ in the above theorem by more general polynomial $P(f)$, where $P(z)$ is defined above?\end{ques}
In this paper, we try to find out the possible answer to the above questions, and also investigated the uniqueness of a zero-order meromorphic function $f$ and its linear q-difference operator sharing sets $S_1$ and $S_2$ with finite weights. The next section includes our main results.
\section{\textbf{Main Results}}
\begin{theo}\label{t4}
	Let $f$ be a transcendental meromorphic of zero order (resp. entire) function, and $\alpha(z)$ be a non-zero small function of $f(z)$. Then for $n>m+3k+4$ (resp. $n >m+k$), $P(f)L_{k}(f,E_q)-a(z)$ has	infinitely many zeros.\end{theo}
\begin{theo}\label{t5}
	Let $f$ be a zero-order transcendental meromorphic (resp. entire) function, and $\alpha(z)$ be a non-zero small function of $f(z)$. Then for $n>m+5k+6$ (resp. $n >m+2k+2$), $P(f)L_{k}(f,\Delta)-\alpha(z)$ has	infinitely many zeros.\end{theo}

\begin{theo}\label{t3}
	Let $f$ be a meromorphic function of zero-order, $q\in\mathbb{C}-\{0\}$, $M(z)=P(f(z))$. If $M(z)$ and $M(qz)$ share a non-zero complex constant $a$ and $\infty$ CM, then $P(f(z))=P(f(qz))$.
\end{theo}

\begin{theo}\label{t1}
	Let $n, k \in\mathbb{N}$, and $f$ be a non-constant zero order meromorphic function such that $E_f(S_1,1)= E_{L_{k}(f,\Delta_{q})}(S_1, 1)$ and $E_f(S_2,0) = E_{L_{k}(f,\Delta_{q})}(S_2, 0)$. If $n\geq 7$, then there exists a constant $t\in \mathbb{C}$ such that $L_{k}(f,\Delta_{q})\equiv tf$, where $t^{n}=1$ and $t\neq -1$.
\end{theo}
\begin{theo}\label{t2}
	Let $n, k \in\mathbb{N}$, and $f$ be a non-constant zero order meromorphic function such that $E_f(S_1,2)= E_{L_{k}(f,\Delta_{q})}(S_1, 2)$ and $E_f(S_2,0) = E_{L_{k}(f,\Delta_{q})}(S_2, 0)$. If $n\geq 6$, then there exists a constant $t\in \mathbb{C}$ such that $L_{k}(f,\Delta_{q})\equiv tf$, where $t^{n}=1$ and $t\neq -1$.\end{theo}
\begin{rem}\label{r1.1}
	Clearly, Theorem \ref{t4} and Theorem \ref{t5} are the improvements of Theorem C and Theorem D, respectively.
\end{rem}
\begin{rem}\label{r1.1a}
	The zero order growth restriction in Theorem \ref{t5} can not be extended to finite order. This can be observed by taking $P(z)=z^{n}$, $f(z)=e^z$, $c=0$ and $q=-n$. Then $P(f(z))[f(qz)-f(z)]-1$ have no zeros.
\end{rem}
\begin{rem}\label{r1.2}
	Obviously, Theorem\ref{t3} is an improvement of Theorem E. \end{rem}
\begin{rem}\label{r1.3}
	If we closely observe the statement of Theorem E, the we see that if a zero order meromorphic function $f(z)$ and its q-shift operator $f(qz)$ share the sets $S_1$ and $S_2$ CM, then $f(z)=tf(qz)$ for a constant $t$ such that $t^n=1$.\end{rem}
The following examples show that Theorem \ref{t1} and \ref{t2} hold for $n=7$ and $n=6$ respectively, for both entire and meromorphic functions.
\begin{exm}
For a positive integer $m$, let $f(z)=z^{m},\; L_k(f,\Delta)=f(qz)-f(z)$, where $q=(1+\omega)^{\frac{1}{m}}$, $\omega^7=1\;(\omega^6=1\; \text{for Theorem \ref{t2}})$. Then one can easily verify that $f(z)$ and $L_{k}(f,\Delta)$ share the sets $S_1$ and $S_2$ and $L_k(f,\Delta)=\omega f(z)$ such that $\omega\neq -1$, $\omega^7=1\; (\omega^6=1\; \text{for Theorem1.5})$.\end{exm}
\begin{exm}
	For a positive integer $m$, let $f(z)=\displaystyle
\frac{1}{z^m},\; L_k(f,\Delta)=f(qz)-f(z)$, where $q=\displaystyle\frac{1}{(1+\omega)^{^\frac{1}{m}}}$, $\omega^7=1\;(\omega^6=1\; \text{for Theorem \ref{t2}})$. Then one can easily verify that $f(z)$ and $L_{k}(f,\Delta)$ share the sets $S_1$ and $S_2$ and $L_k(f,\Delta)=\omega f(z)$ such that $\omega\neq -1$, $\omega^7=1\; (\omega^6=1\; \text{for Theorem1.5})$.\end{exm}

\section{\textbf{Some Lemmas}} We now prove several lemmas which will play key roles in proving the main results of the paper. Let $\mathcal{F}$ and $\mathcal{G}$ be two non-constant meromorphic functions defined by \bea\label{e2.1} F=\left(\frac{f(z)}{a}\right)^{n}\;\; \text{and}\;\; G=\left(\frac{L_{k}(f,\Delta)}{a}\right)^{n}.\eea Henceforth we shall denote by $H$ and $V$ in the following. \be\label{e2.2}H=\left(\frac{\;\;F^{\prime\prime'}}{F}-\frac{2F^{\prime}}{F-1}\right)-\left(\frac{\;\;G^{\prime\prime}}{G^{\prime}}-\frac{2G^{\prime}}{G-1}\right),\ee
\be \label{e2.3} V=\frac{F^{\prime}}{F(F-1)}-\frac{G^{\prime}}{G(G-1)}.\ee

\begin{lem}\label{lem2.1}\cite{Liu & Qi & 2011}
	Let $f$ be a zero-order meromorphic function, and $q(\neq 0),\;c\in\mathbb{C}$. Then \beas m\left(r,\frac{f(qz+c)}{f(z)}.\right)=S(r,f).\eeas \end{lem}
\begin{lem}\label{lem2.2}\cite{Xu & Liu & Cao & 2015}
	Let $f(z)$ be a transcendental meromorphic function of zero order and $q,c$ two non-zero complex constants. Then
	\beas \label{e2.4}T(r,f(qz+c))=T(r,f)+S(r,f),\eeas \beas N(r,\infty;f(qz+c))=N(r,f(z))+S(r,f),\eeas
	\beas N(r,0;f(qz+c))=N(r,f(z))+S(r,f),\eeas
	\beas \ol N(r,\infty;f(qz+c))=N(r,f(z))+S(r,f),\eeas
	\beas \ol N(r,0;f(qz+c))=N(r,f(z))+S(r,f)\eeas on a set of logarithmic density $1$.
\end{lem}

\begin{lem}\label{lem2.3}\cite{7a} If $N\left(r,\displaystyle\frac{1}{f^{(k)}}\mid f\not=0\right)$ denotes the counting function of those zeros of  $f^{(k)}$ which are not the zeros of $f$, where a zero of $f^{(k)}$ is counted according to its multiplicity then $$N\left(r,0;f^{(k)}\mid f\not=0\right)\leq k\ol N(r,\infty;f)+N\left(r,0;f\mid <k\right)+k\ol N\left(r,0f\mid\geq k\right)+S(r,f).$$\end{lem}
\begin{lem}\label{lem2.4}\cite{9} Let $f$ be a non-constant meromorphic function and let \[\mathcal{R}(f)=\frac{\sum\limits _{i=0}^{n} a_{i}f^{i}}{\sum \limits_{j=0}^{m} b_{j}f^{j}}\] be an irreducible rational function in $f$ with constant coefficients $\{a_{i}\}$ and $\{b_{j}\}$ where $a_{n}\not=0$ and $b_{m}\not=0$. Then $$T(r,\mathcal{R}(f))=d\;T(r,f)+S(r,f),$$ where $d=\max\{n,m\}$.\end{lem}
\begin{lem}\label{lem2.5}
	Let $F$ and $G$ be given by (\ref{e2.1}) satisfying $E_F(1,m) = E_G(1,m)$, $0\leq m <\infty$ with $H\not\equiv0$, then \beas N_E^{1)}(r,1;F)\leq N(r,H)+S(r,F)+S(r, G).\eeas Similar inequality holds for $G$ also. \end{lem}
\begin{proof}
	Let $z_{0}$ be a simple 1-point of $F$. Then by a simple calculation, it can be shown that $z_{0}$ is a zero of $H$. Since $m(r,H)=S(r,F)+S(r,G)$, by First fundamental Theorem of Nevalinna, we obtain
	\beas N_E^{1)}(r,1;F)\leq N\left(r,\frac{1}{H}\right)\leq T(r,H)\leq N(r,H)+S(r,F)+S(r,G).\eeas
	\end{proof}
   \begin{lem}\label{lem2.6}
   Let $f$ be a non-constant meromorphic function of zero order, then $S(r,L_{k}(f,\Delta))$ can be replaced by $S(r, f )$.\end{lem}
   \begin{proof} In view of Lemma \ref{lem2.2}, we get \beas T(r,L_{k}(f,\Delta))\leq\sum_{j=0}^{k}T(r,f(q_{j}z))+T(r,f(z))+S(r,f)\leq (k+2)T(r,f)+S(r,f).\eeas Hence, the lemma follows.\end{proof}
   \begin{lem}\label{lem2.7}
   	Let $F$ and $G$, being defined in (\ref{e2.1}) share $(1,m)$, $1\leq m <\infty$ and $(\infty, 0)$. Then
   	\beas \ol N(r,1;F,G)&\leq& \frac{1}{m+1}\big\{ \ol N(r,0;f(z))+\ol N(r,0;L_{k}(f,\Delta))\big\}+\frac{2}{m+1}\ol N(r,\infty; f(z))\\&&+S(r,f).\eeas\end{lem}
   \begin{proof} In view of Lemma \ref{lem2.3} and \ref{lem2.6}, we get \beas \ol N_{*}(r,1;F,G)&=&\ol N_{L}(r,1;F)+\ol N_{L}(r,1;G)\\&\leq& \ol N(r,1;F\mid\geq m+2)+\ol N(r,1;G\mid\geq m+2)\\&\leq&\frac{1}{m+1}\Big\{N(r,0;F^{\prime}\mid F\neq 0)+N(r,0;G^{\prime}\mid G\neq 0)\Big\}\\&\leq&\frac{1}{m+1}\Big\{\ol N(r,0;F)+\ol N(r,0;G)+2\ol N(r,\infty;F)\Big\}+S(r,f)\\&\leq&\frac{1}{m+1}\Big\{ \ol N(r,0;f(z))+\ol N(r,0;L_{k}(f,\Delta))\Big\}+\frac{2}{m+1}\ol N(r,\infty; f(z))\\&&+S(r,f).\eeas\end{proof}

\begin{lem}\label{lem2.9}
	Let $F,$ $G$ be given by $(\ref{e2.1})$ and $V \not\equiv 0$. If $F,$ $G$ share $(1,m)$, and $f,$ $L_{k}(f,\Delta)$ share $(\infty,k),$ where $0\leq m,k<\infty,$ then the poles of $F$ and $G$ are zeros of $V$ and \beas &&(nk+n-1)\ol N(r,\infty;f \mid\geq k+1)\\&\leq& \ol N(r,0;f(z))+\ol N(r,0;L_{k}(f,\Delta))+\ol N_{*}(r,1;F,G)+S(r,f).\eeas
\end{lem}
\begin{proof}
	Since $f(z),$ $L_{k}(f,\Delta)$ share $(\infty;k)$, it follows that $F,$ $G$ share $(\infty;nk)$ and so a pole of $F$ with multiplicity $p(\geq nk + 1)$ is a pole of $G$ with multiplicity $r(\geq nk + 1)$ and vice versa. We note that $F$ and $G$ have no pole of multiplicity $q$ where $nk < q < nk + n$.
	Now using the Milloux theorem [\cite{Hayman & 1964}, p. 55] and Lemma \ref{lem2.6}, we get from the definition of $V$, $ m(r,V)=S(r,f(z)).$
	Therefore, \beas &&(nk+n-1)\ol N(r,\infty;f\mid\geq k+1)\leq\; N(r,0;V)\leq\; T(r,V)+O(1)\\&&\leq\; N(r,\infty;V)+m(r,V)+O(1)\\&&\leq\;\ol N(r,0;F)+\ol N(r,0;G)+\ol N_{*}(r,1;F,G)+S(r,f(z))+S(r,L(z,f))\\&&\leq\; \ol N(r,0;f)+\ol N(r,0;L_{k}(f,\Delta))+\ol N_{*}(r,1;F,G)+S(r,f).\eeas
\end{proof}
\begin{lem}\label{lem2.10}\cite{Yi & 1997, Zhang & 2010}
	If $F$ and $G$ share $(\infty, 0)$ and $V\equiv0$, then $F\equiv G$.\end{lem}
\begin{lem}\label{lem2.11}\cite{Yi & 1999}
	Let $H\equiv0$ and $F$, $G$ share $(\infty, 0)$, then $F$ and $G$ share $(1,\infty)$, $(\infty,\infty)$.\end{lem}
\begin{lem}\label{lem2.13}\cite{Banerjee & 2010}
	Let $F,$ $G$ be two non-constant meromorphic functions sharing $(1,m),$ where $1\leq m<\infty.$ Then \beas && \ol N(r,1;F) + \ol N(r,1;G)- N_{E}^{1)}(r,1;F) +\left(t-\frac{1}{2}\right)\ol N_{*}(r,1;F,G) \\&\leq& \frac{1}{2} [N(r,1;F) + N(r,1;G)].\eeas\end{lem}
\begin{lem}\label{lem2.14}\cite{Lahiri & Banerjee & 2006}
	Suppose $F$, $G$ share $(1,0)$, $(\infty,0)$. If $H\not \equiv 0,$ then, \beas N(r,H) &\leq&\ol  N(r,0;F \mid\geq 2) + \ol N(r,0;G \mid\geq 2)+\ol N_{*}(r,1;F,G) + \ol N_{*}(r,\infty;F,G) \\&&+ \ol N_{0}(r,0;F^{\prime}) + \ol N_{0}(r,0;G^{\prime})+S(r,F)+S(r,G),\eeas where $\ol N_{0}(r,0;F^{\prime})$ is the reduced counting function of those zeros of $F^{\prime}$ which are not the zeros of $F(F-1)$ and $\ol N_{0}(r,0;G^{\prime})$ is similarly defined.\end{lem}
\begin{lem}\label{lem2.15} Let $f$ be a transcendental meromorphic function funtion of zero order and $L_{k}(f,E_q)$ be a linear q-shift polynomial defined in $(\ref{e1.1a})$ and $P(z)$ be defined as in the introduction. Then we have \beas  (n-k-1)T(r,f)+S(r,f)\leq T(r,P(f)L_{k}(f,E_{q}))\leq (n+k+1)T(r,f)+S(r,f.)\eeas
	If $f$ is a transcendental entire function of zero order, then \beas T(r,P(f))L_k(f,E_q)=(n+1)T(r,f)+S(r,f).\eeas\end{lem}
\begin{proof}
	We set a function  $F:=P(f)L_{k}(f, E_q)$, now if $f$ is an entire function of finite order, then \bea && \label{e2.5}T(r,F)= \nonumber m(r,F)\leq m(r,P(f)L_k(f, E_q))+S(r,f)\\&\leq& m(r,P(f)f)+m\left(r,\frac{L_k(f, E_q)}{f}\right)+S(r,f)\leq T(r,P(f)f)+S(r,f).\eea
	On the other hand, using Lemma \ref{lem2.4}, we have \bea \label{e2.6} && (n+1)T(r,f)=\nonumber T(r,P(f)f)+S(r,f)\leq m(r,P(f)f)+S(r,f)\\&\leq & \nonumber m(r,F)+m\left(r,\frac{P(f)f}{F}\right)+S(r,f)\\&\leq&  m(r,F)+m\left(r,\frac{f}{L_k(f,E_q)}\right)+S(r,f)\leq T(r,F)+S(r,f). \eea
	It follows from $(\ref{e2.4})$ and $(\ref{e2.5})$, \beas T(r,F)=(n+1)T(r,f)+S(r,f).\eeas
	If $f$ is a meromorphic function of zero order, then
	\bea T(r,F)=\nonumber T(r,P(f)L_k(f,E_q))&\leq& T(r,P(f))+T(r,L_k(f,E_q))+S(r,f)\\&\leq& (n+k+1)T(r,f)+S(r,f).\eea Also we see that, \bea\label{e2.7} (n+1)T(r,f)+S(r,f)&=&\nonumber T(r,P(f)f)+S(r,f)\\&\leq& \nonumber m(r,P(f)f)+N(r,P(f)f)+S(r,f)\\&\leq&\nonumber m\left(r,F\frac{f}{L_k(f,E_q)}\right)+N\left(r,F\frac{f}{L_k(f,E_q)}\right)+S(r,f)\\&\leq& T(r,F)+(k+2)T(r,f)+S(r,f).\eea
	Equations $(\ref{e2.6})$ and $(\ref{e2.7})$ yield \beas (n-k-1)T(r,f)+S(r,f)\leq T(r,F)\leq (n+k+1)T(r,f)+S(r,f).\eeas\end{proof}
\begin{lem}\label{lem2.16} Let $f$ be a transcendental meromorphic function funtion of zero order and $L_{k}(f,\Delta_q)$ be a linear q-difference polynomial defined in $(\ref{e1.1})$ and $P(z)$ be defined as in the introduction. Then we have \beas T(r,P(f)L_{k}(f,E_{q}))\geq(n-k-1)T(r,f)+S(r,f).\eeas If f is a transcendental entire function of zero order, we have \beas T(r,P(f)L_k(f,\Delta))\geq nT(r,f)+S(r,f).\eeas\end{lem}
\begin{proof} In a similar manner as done in Lemma \ref{lem2.15}, we can prove it. So, we omit the details.\end{proof}
\section{\textbf{Proofs of the theorems}}
\begin{proof}[Proof of Theorem \ref{t4}]
	\textbf{Case 1.} Suppose $f$ be a transcendental meromorphic of zero order  function. Let $P(f)L_{k}(f,E_q)-a(z)$ has	only finitely many zeros.\par By the Second fundamental Theorem of Nevalinna, and Lemmas \ref{lem2.15}, \ref{lem2.2}, we obtain \beas T(r, P(f)L_{k}(f,E_q))&\leq& \ol N(r,\infty;P(f)L_{k}(f,E_q))+\ol N(r,0;P(f)L_{k}(f,E_q))\\&&+\ol N(r,0;P(f)L_{k}(f,E_q)-\alpha(z))+S(r,f)\\&\leq&(k+2)\ol N(r,\infty; f(z))+(m+k+1)T(r,f(z))+S(r,f).\eeas i.e., \beas (n-m-3k-4)T(r,f)\leq S(r,f),\eeas which is a contradiction since $n>m+3k+4$.\par 
	\textbf{Case 2.} Let $f(z)$ is a transcendental entire function. Then using the similar arguments as done in Case 1, we get \beas (n-m-k)T(r,f)\leq S(r,f),\eeas which contradicts $n>m+k$. Hence, the theorem is proved.
\end{proof}
\begin{proof}[Proof of Theorem \ref{t5}]
	Using Lemma \ref{lem2.16} and \ref{lem2.2}, the proof of this theorem can be done in the line of proof of the Theorem \ref{t4}. So we the details.
\end{proof}
\begin{proof}[Proof of Theorem \ref{t3}]
	Let $G_*(z)=\displaystyle\frac{M(z)}{a}$, then we know $G_*(z)$ and $G_*(qz)$ share $1$ and $\infty$ CM. Since the order of $f$ is zero, it follows that
	\beas \frac{G_*(z)-1}{G_*(qz)-1}=h,\;\; \text{where h is a non-zero constant}.\eeas Rewriting the above equation, we obtain \bea \label{e3.3} G_*(z)+\frac{1}{h}-1=\frac{G_*(qz)}{h}.\eea Suppose $h\neq 1$. Keeping in view of \ref{e2.4} and \ref{e3.3}, the Second Fundamental theorem of Nevalinna yields \beas nT(r,f(z))&=&T(r, G_*(z))\leq \ol N(r,\infty;G_*)+\ol N(r, 0;G_*)+\ol N(r, 1-\frac{1}{h};G_*)+S(r,f)\\&\leq& N(r,\infty;f(z))+\ol N(r, 0;P(f))+\ol N(r, 0;G_*(qz))+S(r,f)\\&\leq& (2m+1)T(r,f)+S(r,f),\eeas which contradicts $n>2m+1$. Hence $h=1$, which implies that $G_{*}(z)=G_*(qz)$, i.e., $P(f(z))=tP(f(qz))$ for a constant $t$ satisfying $t^n=1$. 
\end{proof}

\begin{proof}[Proof of Theorem \ref{t1}]
Let $F$ and $G$ be two functions defined in (\ref{e2.1}).\par  Since $E_{f}(S_1,1)=E_{L_{k}(f,\Delta)}(S_1,1)$ and $E_{f}(S_2,0)=E_{L_{k}(f,\Delta)}(S_2,0)$, it follows that $F$, $G$ share $(1,1)$ and $(\infty,0).$\par
We now consider the following two cases.\par 
\noindent {\bf{Case 1:}} Suppose $H\not\equiv 0$. Then $F\not\equiv G.$ So, by Lemma \ref{lem2.10}, it follows that $V\not\equiv 0.$\par
By the Second Fundamental Theorem of Nevalinna, we have
\bea \label{e3.1}&&T(r,F)+T(r,G)\nonumber \\&\leq& \ol N(r,1;F)+\ol N(r,0;F)+\ol N(r,\infty;F)+\ol N(r,1;G)+\ol N(r,0;G)+\ol N(r,\infty;G)\nonumber\\&&-\ol N_0(r,0;F^{\prime})-\ol N_0(r,0;G^{\prime})+S(r,F)+S(r,G).\eea
Applying Lemmas \ref{lem2.5}, \ref{lem2.13} and \ref{lem2.14} to the \ref{e3.1}, we obtain, \beas && \frac{n}{2}\Big(T(r,F)+T(r,G)\Big)\\&\leq& N_2(r,0;F)+N_2(r,0;G)+3\ol N(r,\infty; F)-\left(m-\frac{3}{2}\right)\ol N_*(r,1;F,G)\\&&+S(r,F)+S(r,G)\\&\leq& 2\{\ol N(r,0;F)+\ol N(r,0;G)\}+3\ol N(r,\infty;F)-\left(m-\frac{3}{2}\right)\ol N_{*}(r,1;F,G)\\&&+S(r,F)+S(r,G).\eeas
Using Lemma \ref{lem2.7} with $m=1$ and Lemma \ref{lem2.6}, the above inequality becomes
\beas  \frac{n}{2}\Big(T(r,F)+T(r,G)\Big)\leq \frac{9}{4}\Big(\ol N(r,0;F)+\ol N(r,0;G)\Big)+\frac{7}{2}\ol N(r,\infty;F)+S(r,f).\eeas
Again applying Lemmas \ref{lem2.4}, \ref{lem2.6}, \ref{lem2.7} and \ref{lem2.9} with $m=1$, $k=0$, the above equation turns into 
\beas \left(\frac{n}{2}-\frac{21}{4(n-2)}-\frac{9}{4}\right)\left(T(r,f)+T(r,L_{k}(f,\Delta))\right)\leq S(r,f),\eeas which contradicts $n\geq 7$.\par 
{\bf{Case 2:}} Suppose $H\equiv 0$. After integration we get, \be\label{e3.2} F\equiv \frac{\alpha G+\beta}{\gamma G+\delta},\ee where $\alpha,\beta,\gamma,\delta$ are complex constants satisfying $\alpha\delta-\beta\gamma\neq 0.$\par As $F$, $G$ share $(\infty,0)$, by Lemma \ref{lem2.11}, it follows that $f(z)$, $L_{k}(f,\Delta)$ share $(1,\infty)$ and $(\infty, \infty)$.\par 
{\bf{Subcase 2.1:}} Let $\alpha\gamma \neq 0$. Then $F-\displaystyle\frac{\alpha}{\gamma}=\frac{-(\alpha\delta-\beta\gamma)}{\gamma(\gamma G+\delta)}\neq 0.$ \par
Therefore, by the Second Fundamental Theorem of Nevallina, we get \beas nT(r,f)&\leq& \ol N(r,0;F)+\ol N(r,\infty;F)+\ol N(r,\frac{\alpha}{\gamma};F)+S(r,F)\\ &\leq & 2T(r,f)+S(r,f)\eeas which is a contradiction since $n\geq 7$.\par
{\bf{Subcase 2.2:}} Suppose that $\alpha\gamma=0$. Since $\alpha\delta-\beta\gamma\neq 0,$ both $\alpha$ and $\gamma$ are not zero simultaneously.\par
{\bf{Subcase 2.2.1:}} Suppose $ \alpha\neq 0$ and $\gamma=0.$ Then (\ref{e3.2}) becomes $ F \equiv A G+B $, where $ A=\displaystyle\frac{\alpha}{\delta}$ and $B=\displaystyle\frac{\beta}{\delta}.$\par
{\bf{Subcase 2.2.1.1:}} Let $F$ has no $1$-point. Then by the Second Fundamental Theorem, we get \beas T(r,F)\leq \ol N(r,0;F)+\ol N(r,1;F)+\ol N(r,\infty;F)+S(r,F)\eeas $or,$ \beas (n-2)T(r,f)\leq S(r,f),\eeas which is a contradiction.\par
{\bf{Subcase 2.2.1.2:}} Let $F$ has some $1$-point. Then $A+B=1$.\par
{\bf{Subcase 2.2.1.2.1:}} Suppose $A\neq1$. Then $ F \equiv A G+1-A$.\par
Therefore, by the Second Fundamental Theorem, we get \beas T(r,F)&\leq& \ol N(r,0;F)+\ol N(r,\infty;F)+\ol N(r,1-\alpha;F)+S(r,F)\\&\leq& \ol N(r,0;F)+\ol N(r,\infty;F)+\ol N(r,0;G)++S(r,F)\\&\leq& 3T(r,f)+S(r,f).\eeas $i.e.,$ \beas (n-3)T(r,f)\leq S(r,f),\eeas which is again a contradiction since $n\geq 7$.\par
{\bf{Subcase 2.2.1.2.2:}} Suppose $A=1$. Then $F\equiv G.$ Thus we have $L_{k}(f,\Delta)\equiv tf(z)$, where $t^{n}=1$ and $t\neq -1$\par
{\bf{Subcase 2.2.2:}} Suppose $\alpha=0$ and $\gamma\neq 0$.\par

Then (\ref{e3.2}) becomes \beas F\equiv \frac{1}{C G+D},\;\;\text{where}\;\; C=\frac{\gamma}{\beta}, \;\; D=\frac{\delta}{\beta}.\eeas 
{\bf{Subcase 2.2.2.1:}} Let $F$ has no $1$-point. Then applying the second fundamental theorem to $F$, we have \beas nT(r,f)&\leq &\ol N(r,\infty;F)+\ol N(r,0;F)+\ol N(r,1;F)+S(r,F)\\&\leq & 2T(r,f)+ S(r,f),\eeas which is a contradiction.\par
{\bf{Subcase 2.2.2.2:}} Suppose that $F$ has some $1$-point. Then $C+D=1$.\par
{\bf{Subcase 2.2.2.2.1:}} Suppose $C=1$. Then $D=0$ and thus $FG\equiv 1$. $i.e.$, $f(z)L_{k}(f,\Delta)\equiv ta^{2}$, where $t^{n}=1$. Since $F$ and $G$ share $(\infty,\infty)$, so we have $N\left(r, \displaystyle\frac{L_{k}(f,\Delta)}{f}\right)= N(r, 0; f)$ and so in view of Lemma \ref{lem2.1} and \ref{lem2.6}, we have \beas && 2T(r,f)+S(r,f)\leq T\left(r,\frac{ta^{2}}{f^{2}}\right)\leq T\left(r,\frac{L_{k}(f,\Delta)}{f}\right)+S(r,f)\\&\leq& N\left(r,\frac{L_{k}(f,\Delta)}{f}\right)+S(r,f)\leq N(r,0;f)+S(r,f)\leq T(r,f)+S(r,f), \eeas which is a contradiction.
\end{proof}

\begin{proof}[Proof of Theorem \ref{t2}]
	The proof of this theorem can be done in a similar manner as done in Theorem \ref{t1}, so we omit the details.
\end{proof}
\textbf{Acknowledgment.}\\
The authors did not receive support from any organization for the submitted work.

\end{document}